\documentclass[11pt]{article}

\usepackage{amssymb}
\usepackage{amsmath}
\usepackage{amsthm}
\usepackage{mathrsfs}
\usepackage{graphicx}
\usepackage[utf8]{inputenc}
\usepackage{multirow}
\usepackage{floatrow}
\usepackage{mathrsfs}

\newtheorem{thm}{Theorem}[section]

\newtheorem{cor}[thm]{Corollary}

\newenvironment{Pf}{\noindent{\bf Proof. }}{\hfill $\blacksquare$ \\}

\def\ov{\overline}
\def\lf{\lfloor}
\def\rf{\rfloor}
\def\lc{\lceil}
\def\rc{\rceil}

\def\ts{\textstyle}
\def\pr{\prime}
\def\sm{\setminus}
\def\b{\beta}
\def\t{\tilde}
\def\C{\mathcal C}
\def\D{\mathscr D}
\def\De{\Delta}
\def\K{\mathcal K}

\def\T{\mathcal T}
\def\W{\mathcal W}
\def\F{\mathcal F}
\def\P{\mathcal P}

\begin{document}

\title{\bf Measures of closeness to cordiality for graphs}
\author{
{\bf Anand Brahmbhatt}\thanks{Department of Computer Science, Princeton University, 35 Olden Street, Princeton, NJ 08540-5233, USA. \newline {\tt e-mail:ab7728@princeton.edu}} \qquad
{\bf Kartikeya Rai}\thanks{Centre for Mathematical Sciences, Wilberforce Road, Cambridge, CB3 0WA, United Kingdom. \newline {\tt e-mail:kr547@cantab.ac.uk}} \qquad
{\bf Amitabha Tripathi}\thanks{Department of Mathematics, Indian Institute of Technology, Hauz Khas, New Delhi -- 110016, India. \newline {\tt e-mail:atripath@maths.iitd.ac.in}}\:\:\thanks{\it Corresponding author}
}
\date{}
\maketitle

\begin{abstract}
\noindent A graph $G$ is cordial if there exists a function $f$ from the vertices of $G$ to $\{0,1\}$ such that the number of vertices labelled $0$ and the number of vertices labelled $1$ differ by at most $1$, and if we assign to each edge $xy$ the label $|f(x)-f(y)|$, the number of edges labelled $0$ and the number of edges labelled $1$ also differ at most by $1$. We introduce two measures of how close a graph is to being cordial, and compute these measures for a variety of classes of graphs. 
\end{abstract}


\vskip 20pt

\section{Introduction} \label{intro}
\vskip 10pt

Graph labellings were introduced by Rosa \cite{Ros67} in a bid to attack the conjecture of Ringel \cite{Rin64} that ${\K}_{2n+1}$ can be decomposed into $2n+1$ subgraphs that are all isomorphic to a given tree with $n$ edges. Rosa called a function $f$ a $\b$-{\it valuation\/} of a graph $G$ with $m$ edges if $f$ is an injection from the vertices of $G$ to the set $\{0,1,\ldots,m\}$ such that, when each edge $xy$ is assigned the label $|f(x)-f(y)|$, the resulting edge labels are distinct. Golomb \cite{Gol72} subsequently called such labellings {\it graceful\/}. Ringel's conjecture has been proved for all sufficiently large $n$ recently by Keevash and Staden \cite{KS20} in April 2020, and independently by Montgomery, Pokrovskiy and Sudakov \cite{MPS21} in January 2021. The conjecture of Ringel followed by the Graceful Tree Conjecture has spurred a great deal of activity in the area of Graph Labellings over the past six decades; for an updated and dynamic survey, see Gallian \cite{Gal24}. Among the numerous offshoots of the graceful labelling of a graph, and among the most prominent is the cordial labelling of a graph, introduced by Cahit \cite{Cah87} in 1987.  

Let $G$ be a finite, simple graph. Let $f: V(G) \to \{0,1\}$, and let $\ov{f}: E(G) \to \{0,1\}$ be the induced mapping given by 
\[  \ov{f}(xy) = \left| f(x)-f(y) \right|. \]
We say that $f$ is a {\it cordial\/} labelling of $G$ if the number of vertices labelled $0$ and the number of vertices labelled $1$ differ by at most $1$, and the number of edges labelled $0$ and the number of edges labelled $1$ differ at most by $1$. A graph is said to be {\it cordial\/} if it admits a cordial labelling. The vastness of literature relating to cordial graphs may be measured by the contents in \cite[pp. 89--107]{Gal24}. 

For $i \in \{0,1\}$, let $v_i(f)$ denote the number of vertices labelled $i$ and let $e_i(f)$ denote the number of edges labelled $i$. Let 
\[ {\De}_v(f) = \left| v_0(f)-v_1(f) \right| \;\;\text{and}\;\; {\De}_e(f) = \left| e_0(f)-e_1(f) \right|. \] 
Hence, $f$ is a cordial labelling of $G$ precisely when ${\De}_v(f) \le 1$ and ${\De}_e(f) \le 1$. 

Chartrand, Lee and Zhang \cite{CLZ06} introduced the notion of uniform cordiality. They called a labelling $f: V(G) \to \{0,1\}$ {\it friendly\/} if ${\De}_v(f) \le 1$. A graph $G$ for which every friendly labelling is cordial is called {\it uniformly cordial\/}. They proved that a connected graph of order $n \ge 2$ is uniformly cordial if and only if $G={\K}_3$ or $G={\K}_{1,n-1}$, $n$ even.

Riskin \cite{Ris07a,Ris07b} introduced two measures of the noncordiality of a graph and determined these measures for certain classes of graphs. He defined the {\it cordial vertex deficiency\/} of a graph $G$ (denoted by $\text{cvd}(G)$) as the minimum number of vertices, taken over all labellings of $G$ for which ${\De}_e(f) \le 1$, which needs to be added to $G$ such that the resulting graph is cordial. Analogously, he defined the {\it cordial edge deficiency\/} (denoted by $\text{ced}(G)$) of a graph $G$ as the minimum number of edges, taken over all labellings of $G$ for which ${\De}_v(f) \le 1$, which needs to be added to $G$ such that the resulting graph is cordial. 

In this paper, we define two measures of closeness to cordiality of a graph as follows: 
\begin{eqnarray}
{\D}_1(G) = \min_{f} \Big\{ {\De}_v(f) + {\De}_e(f) \Big\}. \\
{\D}_2(G) = \min_{ {\De}_v(f) \le 1} \Big\{ {\De}_e(f) \Big\}. 
\end{eqnarray}
The first measure is over all labellings $f: V(G) \to \{0,1\}$, whereas the second measure is over those labellings $f: V(G) \to \{0,1\}$ for which ${\De}_v(f) \le 1$. So if ${\D}_1(G) \le 1$ or if ${\D}_2(G) \le 1$, then $G$ is cordial. Thus, ${\D}_1$ and ${\D}_2$ may both be viewed as measures of cordiality. We note that the measure ${\D}_2(G)$ is closely related to the measure $\text{ced}(G)$.   

This paper is divided into sections as follows. We provide sharp upper bounds for the join of two graphs with respect to each of the two measures in Section \ref{graph_join}. In subsequent sections, we study these measures for several well known classes of graphs, some of which are the join of two well known classes of graphs. More specifically, we provide exact values for both measures, for trees (Section \ref{tree}), complete graphs (Section \ref{complete}), complete $r$-partite graphs (Section \ref{r-partite}), cycles (Section \ref{cycle}), wheels (Section \ref{wheel}), and fans (Section \ref{fan}), except that we only provide bounds for ${\D}_1(G)$ in the case of complete $r$-partite graphs. A summary of the results is given in Table \ref{results}.  We close the paper with a few directions of enquiry (Section \ref{CR}). 

\begin{table}[H]
\centering
\renewcommand{\arraystretch}{1.1}
\begin{tabular}{|c|c|} \hline
$G$ & ${\D}_1\big(G\big)$ \\ \hline \hline
$G_1+G_2$ & $\le {\D}_1(G_1)\,{\D}_1(G_2)+{\D}_1(G_1)+{\D}_1(G_2)$ \\ \hline
${\T}_n$ & $1$ \\ \hline 
${\K}_n$ & $\begin{cases} a+\frac{1}{2}(n-a^2) & \:\text{if}\: n \in \{a^2+2t: 0 \le t \le a\}, \\ 2a-1 & \:\text{if}\:  n=a^2+1, \\ a+1+\frac{1}{2}\big((a+1)^2-n\big) & \:\text{if}\: n \in \{a^2+2t+1: 1 \le t \le a-1\}.  \end{cases}$ \\ \hline 
& \multicolumn{1}{l|}{If $s$ among $n_1,\ldots,n_r$ are odd, then} \\   
${\K}_{n_1,\ldots,n_r}$ & $\begin{cases} = \sqrt{s} & \text{ if } s \text{ is a square}, \\ \le \frac{1}{2} \big(s+1 - (2a-1)^2\big) & \text{ if } (2a)^2 < s < (2a+1)^2 \text{ and } s \text{ is even}, \\
\le \frac{1}{2} \big(s+1 - (2a-2)^2\big) & \text{ if } (2a)^2 < s < (2a+1)^2 \text{ and } s \text{ is odd}, \\ \le \frac{1}{2} \big(s+1 - (2a-1)^2\big) & \text{ if } (2a+1)^2 < s < (2a+2)^2 \text{ and } s \text{ is even}, \\
\le \frac{1}{2} \big(s+1 - (2a)^2\big) & \text{ if } (2a+1)^2 < s < (2a+2)^2 \text{ and } s \text{ is odd}. \end{cases}$ \\ 
&  \multicolumn{1}{l|}{$\in \left[ \:\lf \sqrt{s} \rf, 3 \lf \sqrt{s} \rf \:\right]$} \\ \hline 
${\C}_n$ & $\begin{cases} 0 & \:\text{if}\: n \equiv \:0\!\!\pmod{4}, \\ 2 & \:\text{if}\: n \not\equiv \:0\!\!\pmod{4}. \end{cases}$ \\ \hline
${\W}_n$ & $\begin{cases} 0 & \:\text{if}\: n \equiv \:2\!\!\pmod{4}, \\ 1 & \:\text{if}\: n \equiv \:1\!\!\pmod{2}, \\ 2 & \:\text{if}\: n \equiv \:0\!\!\pmod{4}. \end{cases}$  \\ \hline
${\F}_{m,n}$ & $\begin{cases} 2 & \:\text{if}\: m \:\text{is odd}, n \:\text{is even}, \\ 1 & \:\text{otherwise}. \end{cases}$ \\ \hline
\end{tabular}
\caption{Summary of results on ${\D}_1(G)$} \label{results}
\end{table}
\vskip 10pt

\begin{table}[H]
\centering
\renewcommand{\arraystretch}{1.4}
\begin{tabular}{|c|c|} \hline
$G$ & ${\D}_2\big(G\big)$ \\ \hline \hline
$G_1+G_2$ & $\le {\D}_2(G_1)+{\D}_2(G_2)+1$ \\ \hline
${\T}_n$ & $1-(n \bmod{2})$ \\ \hline 
${\K}_n$ &  $\left\lf \tfrac{n}{2} \right\rf$ \\ \hline 
${\K}_{n_1,\ldots,n_r}$ & $\left\lf \tfrac{s}{2} \right\rf$ \\ \hline  
${\C}_n$ &  $\begin{cases} 0 & \:\text{if}\: n \equiv \:0\!\!\pmod{4}, \\ 1 & \:\text{if}\: n \equiv \:1, 3\!\!\pmod{4}, \\ 2 & \:\text{if}\: n \equiv \:2\!\!\pmod{4}. \end{cases}$  \\ \hline
${\W}_n$ & $\begin{cases} 0 & \:\text{if}\: n \not\equiv \:0\!\!\pmod{4}, \\ 2 & \:\text{if}\: n \equiv \:0\!\!\pmod{4}. \end{cases}$ \\ \hline
${\F}_{m,n}$ & $\begin{cases} 0 & \:\text{if}\: m \:\text{is even}, n \:\text{is odd}, \\ 1 & \:\text{otherwise}. \end{cases}$ \\ \hline
\end{tabular}
\caption{Summary of results on ${\D}_2(G)$} \label{results_2}
\end{table}
\vskip 10pt

\section{Bound on $\D_1(G_1+G_2)$ and $\D_2(G_1+G_2)$} \label{graph_join}
\vskip 10pt

By the join $G_1+G_2$ of graphs $G_1$ and $G_2$ with disjoint vertex sets $V_1=V(G_1)$ and $V_2=V(G_2)$, we mean the graph $G$ with vertex set $V(G)=V_1 \cup V_2$ and edge set $E(G)$ consisting of edges from $G_1$ and $G_2$ together with edges with one endpoint from $V_1$ and the other from $V_2$. For instance, the star graph ${\K}_{1,n}=\ov{{\K}_n}+{\K}_1$, the complete $r$-partite graph ${\K}_{n_1,\ldots,n_r}=\ov{{\K}_{n_1}}+\cdots+\ov{{\K}_{n_r}}$, the wheel graph ${\W}_n={\C}_{n-1}+{\K}_1$, and the fan graph ${\F}_{m,n}=\ov{{\K}_m}+{\P}_n$. In this section, we find upper bounds for ${\D}_1(G_1+G_2)$ and for ${\D}_2(G_1+G_2)$ in terms of ${\D}_1(G_1)$, ${\D}_1(G_2)$, ${\D}_2(G_1)$ and ${\D}_2(G_2)$. 
\vskip 5pt

\begin{thm} \label{join}
For graphs $G_1$ and $G_2$ with disjoint vertex sets, 
\begin{eqnarray*}
& \D_1(G_1+G_2) \le  \D_1(G_1)\D_1(G_2) + \D_1(G_1) + \D_1(G_2); \\
& \D_2(G_1+G_2) \le \D_2(G_1) + \D_2(G_2) + 1.
\end{eqnarray*}
\end{thm}

\begin{Pf}
There is a one-to-one correspondence between labellings $f: V(G_1+G_2) \to \{0,1\}$ and pairs of labelling $f_1: V(G_1) \to \{0,1\}$ and $f_2: V(G_2) \to \{0,1\}$. The following hold by the triangle inequality.  

\noindent We have 
\begin{equation} \label{vertex}
{\De}_v(f) = | v_0(f) - v_1(f) | \le | v_0(f_1) - v_1(f_1) | + | v_0(f_2) - v_1(f_2) | \le {\De}_v(f_1) + {\De}_v(f_2). 
\end{equation}

\noindent Since $e_0(f) = e_0(f_1)+e_0(f_2)+v_0(f_1) v_0(f_2)+v_1(f_1) v_1(f_2)$ and $e_1(f) = e_1(f_1)+e_1(f_2)+v_0(f_1) v_1(f_2)+v_1(f_1) v_0(f_2)$, we have 
\begin{eqnarray}
{\De}_e(f) & = & | e_0(f) - e_1(f) |  \nonumber \\ 
& \le & | e_0(f_1) - e_1(f_1) | + | e_0(f_2) - e_1(f_2) | + \left|\big(v_0(f_1) - v_1(f_1)\big) \big(v_0(f_2) - v_1(f_2)\big)\right| \nonumber \\
& \le & {\De}_v(f_1)\,{\De}_v(f_2) + {\De}_e(f_1) + {\De}_e(f_2). \label{edge}
\end{eqnarray}
\vskip 5pt

\noindent We use the above inequalites to find upper bounds for $\D_1(G_1+G_2)$ and $\D_2(G_1+G_2)$. If $\D_1(G_1)={\De}_v(f_1^{\star}) + {\De}_e(f_1^{\star})$, $\D_1(G_2)={\De}_v(f_2^{\star}) + {\De}_e(f_2^{\star})$ and $f^{\star}$ is the labelling of $G_1+G_2$ corresponding to the pair of labellings $f_1^{\star}, f_2^{\star}$ of $G_1,G_2$ respectively, then by eqn.~\eqref{vertex} and eqn.~\eqref{edge} we have 
\begin{eqnarray*}
\D_1(G_1+G_2) & = & \min_{f_1,\,f_2} \big\{ {\De}_v(f) + {\De}_e(f) \big\} \\
& \le & {\De}_v(f^{\star}) + {\De}_e(f^{\star}) \\
& \le & {\De}_v(f_1^{\star}) + {\De}_v(f_2^{\star}) + {\De}_e(f_1^{\star}) + {\De}_e(f_2^{\star}) + {\De}_v(f_1^{\star})\,{\De}_v(f_2^{\star}) \\
& \le & \D_1(G_1) + \D_1(G_2) + \D_1(G_1)\D_1(G_2).
\end{eqnarray*}
\vskip 5pt

\noindent To obtain an upper bound for $\D_2(G_1+G_2)$, we need to consider cases. Note that 
\begin{eqnarray*}
\D_2(G_1+G_2) = \min_{{\De}_v(f) \le 1} \big\{ {\De}_e(f) \big\}. 
\end{eqnarray*}
\vskip 5pt

\noindent Suppose both $|V(G_1)|$ and $|V(G_2)|$ are odd.  Thus, from eqn.~\eqref{edge}, 
\begin{eqnarray*}
\D_2(G_1+G_2) & = & \min_{{\De}_v(f)=0} \big\{ {\De}_e(f) \big\} \\
& \le & \min_{\substack{v_0(f_1)-v_1(f_1)=1 \\ v_0(f_2)-v_1(f_2)=-1}} \big\{ {\De}_e(f) \big\} \\
& \le & \min_{\substack{v_0(f_1)-v_1(f_1)=1 \\ v_0(f_2)-v_1(f_2)=-1}} \big\{ {\De}_e(f_1) + {\De}_e(f_2) + 1 \big\} \\
& = & \min_{v_0(f_1)-v_1(f_1)=1} \big\{ {\De}_e(f_1) \big\} + \min_{v_0(f_2)-v_1(f_2)=-1} \big\{ {\De}_e(f_2) \big\} + 1 \\
& = & \D_2(G_1) + \D_2(G_2) + 1,
\end{eqnarray*} 
where the last equality follows by inverting the labels of vertices of $G_1$ and $G_2$.
\vskip 5pt

\noindent If both $|V(G_1)|$ and $|V(G_2)|$ are even, then 
\begin{eqnarray*}
\D_2(G_1+G_2) & = & \min_{{\De}_v(f)=0} \big\{ {\De}_e(f) \big\} \\
& \le & \min_{\substack{{\De}_v(f_1)=0 \\ {\De}_v(f_2)=0}} \big\{ {\De}_e(f) \big\} \\
& \le & \min_{\substack{{\De}_v(f_1)=0 \\ {\De}_v(f_2)=0}
} \big\{ {\De}_e(f_1) + {\De}_e(f_2) \big\} \\
& \le & \D_2(G_1) + \D_2(G_2).
\end{eqnarray*} 
\vskip 5pt

\noindent If $|V(G_1)|$, $|V(G_2)|$ are of opposite parity, say with $|V(G_1)|$ even, then 
\begin{eqnarray*}
\D_2(G_1+G_2) & = & \min_{{\De}_v(f)=1} \big\{ {\De}_e(f) \big\} \\
& \le & \min_{\substack{{\De}_v(f_1)=0 \\ {\De}_v(f_2)=1}} \big\{ {\De}_e(f) \big\} \\
& \le & \min_{\substack{{\De}_v(f_1)=0 \\ {\De}_v(f_2)=1}} \big\{ {\De}_e(f_1) + {\De}_e(f_2) \big\} \\
& \le & \D_2(G_1) + \D_2(G_2).
\end{eqnarray*} 
\end{Pf}
\vskip 10pt

\section{The Trees ${\T}_n$} \label{tree}
\vskip 10pt

Cahit \cite{Cah87} proved that every tree is cordial. If ${\T}_n$ denotes a tree of order $n$, this implies ${\D}_1({\T}_n) \le 2$ and ${\D}_2({\T}_n) \le 1$. We give a direct proof relying on the fact that every non-trivial tree has at least two leaves (for instance, the endpoints of every maximal path in the tree), and the removal of each leaf from a tree results in a tree. 
\vskip 5pt

\begin{thm} Let ${\T}_n$ denote a tree with $n$ vertices. Then 
\[ \D_1\big({\T}_n\big) = 1 \quad \text{and} \quad \D_2\big({\T}_n\big) = 1 - (n \bmod{2}). \]
\end{thm}

\begin{Pf}
Let ${\T}_n$ be a tree of order $n$, so that its size is $n-1$. We exhibit a labelling $f$ of ${\T}_n$ such that
\[ v_0(f) - v_1(f) = n \bmod{2} \quad \text{and} \quad e_0(f) - e_1(f) = \begin{cases}
                                                                                                                 \pm 1 & \text{ if } n \equiv 0 \!\!\!\!\pmod{2}, \\
                                                                                                                  0 & \text{ if } n \equiv 1 \!\!\!\!\pmod{2}.
                                                                                                               \end{cases}
\]
Note that such a labelling attains the values of $\D_1\big({\T}_n\big)$ and $\D_2\big({\T}_n\big)$ in the statement of the theorem. That there can be no labelling with smaller values follows from parity arguments.
\vskip 5pt

\noindent We prove the existence of such a labelling by induction on $n$. As the base cases, for $n=1$ choose $f$ to be the labelling which labels the vertex with $0$, and for $n=2$ choose $f$ to be the labelling which labels one vertex with $0$ and the other with $1$.
\vskip 5pt
   
\noindent We assume that each tree of order less than $n$ has a labelling satisfying the above mentioned conditions. Let ${\T}_n$ be a tree with $n \ge 3$, and let $x, y$ be leaves in ${\T}_n$. Then ${\T}_n \sm \{x,y\}$ is a tree of order $n-2$. Let $f$ be a labelling of ${\T}_n \sm \{x,y\}$ satisfying the conditions above. We extend this labelling to ${\T}_n$ by labelling $x$ and $y$ differently; if $\ov{f}$ denotes any such extension to ${\T}_n$, then $v_0(\ov{f})-v_1(\ov{f})=v_0(f)-v_1(f)$. 
\vskip 5pt

Labelling $x$ and $y$ is done according to the following procedure. Since $x, y$ are leaves in ${\T}_n$, each has a unique neighbour, say $x^{\pr}, y^{\pr}$, respectively. If $x^{\pr}, y^{\pr}$ have the same label in ${\T}_n \sm \{x, y\}$, then the additional edges $xx^{\pr}, yy^{\pr}$ have different labels, so that $e_0(\ov{f})-e_1(\ov{f})=e_0(f)-e_1(f)$. If $x^{\pr}, y^{\pr}$ have different labels in ${\T}_n \sm \{x, y\}$ and $n$ is even, then label $x, x^{\pr}$ the same (and $y, y^{\pr}$ the same) if $e_0(f)-e_1(f)=-1$ and $x, x^{\pr}$ different (and $y, y^{\pr}$ different) if $e_0(f)-e_1(f)=1$. Thus, the difference $e_0-e_1$ alternates between $-1$ and $+1$ as $n$ runs through the positive even integers, starting with $-1$ for $n=2$.  
\vskip 5pt

The case where $x^{\pr}, y^{\pr}$ have different labels in ${\T}_n \sm \{x, y\}$ and $n$ is odd remains to be resolved. Let $n$ be odd, and let us assume the existence of the said labelling for trees of all orders $<n$. Thus, there exists a labelling $\t{f}$ on a tree ${\T}_{n-2}$ for which $v_0(\t{f})-v_1(\t{f})=1$ and $e_0(\t{f})-e_1(\t{f})=0$. Adding a leaf $x$ to ${\T}_{n-2}$ and labelling $x$ as $1$ results in a tree ${\T}_{n-1}$ with labelling $f^{\pr}$, with ${\De}_v(f^{\pr})=0$ and ${\De}_e(f^{\pr})=1$. Further adding a leaf $y$ to ${\T}_{n-1}$ results in a tree ${\T}_n$. If $\ov{y}$ is the vertex adjacent to $y$ in ${\T}_n$, we may label $y$ with the same label as $\ov{y}$ has in ${\T}_{n-1}$ if $e_0(f^{\pr})-e_1(f^{\pr})=-1$ and the opposite label if $e_0(f^{\pr})-e_1(f^{\pr})=1$. The extended labelling $f$ on ${\T}_n$ has the desired property.  
\vskip 5pt

\noindent This completes the proof of the claim by induction, and thus completes the proof of the theorem.
\end{Pf}
\vskip 10pt

\section{The complete graphs ${\K}_n$} \label{complete}
\vskip 10pt

Cahit \cite{Cah87} proved that the complete graph ${\K}_n$ is cordial if and only if $n \le 3$. We determine both ${\D}_1$ and ${\D}_2$ for complete graphs ${\K}_n$. 
\vskip 5pt

\begin{thm} \label{K_1}
If $a^2 \le n < (a+1)^2$, then
\begin{itemize}
\item[{\rm (i)}] 
\[ 
{\D}_1\big({\K}_n\big) = 
\begin{cases} 
a+\frac{1}{2}(n-a^2) & \:\:\mbox{if}\:\: n \in \{a^2+2t: 0 \le t \le a\}, \\
2a-1 & \:\:\mbox{if}\:\: n=a^2+1, \\
a+1+\frac{1}{2}\big((a+1)^2-n\big) & \:\:\mbox{if}\:\: n \in \{a^2+2t+1: 1 \le t \le a-1\}. \\
\end{cases} 
\]
\item[{\rm (ii)}]
\[ {\D}_2\big({\K}_n\big) = \left\lf \tfrac{n}{2} \right\rf. \] 
\end{itemize}
\end{thm}

\begin{Pf}
Let $f: V\big({\K}_n\big) \to \{0,1\}$, and let $v_0(f)=k$. By interchanging the labels $0$ and $1$, we may assume $0 \le k \le \lf\frac{n}{2}\rf$. Thus,  
\begin{eqnarray*} \label{K_n_basic}
v_0(f) = k, \;\; v_1(f) = n-k, \quad {\De}_v(f) = n-2k. \\
e_0(f) = \ts{k \choose 2} + \ts{n-k \choose 2}, \;\; e_1(f) = k(n-k), \quad {\De}_e(f) = \left| \ts{k \choose 2} + \ts{n-k \choose 2} - k(n-k) \right|. 
\end{eqnarray*}
The expression for ${\De}_e(f)$ can be written as $\frac{1}{2}\left| (n-2k)^2-n \right|$. Therefore, 
\begin{eqnarray} \label{K_n_D1,2basic}
{\D}_1\big({\K}_n\big) = \min \Big\{ (n-2k) + \tfrac{1}{2}\left| (n-2k)^2-n \right| : 0 \le k \le \left\lf\tfrac{n}{2}\right\rf \Big\}, \\
{\D}_2\big({\K}_n\big) = \min \Big\{ \tfrac{1}{2}\left| (n-2k)^2-n \right| : | n-2k | \le 1 \Big\} = \left\lf \tfrac{n}{2} \right\rf. 
\end{eqnarray}

\noindent {\sc Case} I. If $n-2k \ge \sqrt{n}$, then the expression in eqn.~\eqref{K_n_D1,2basic} is $(n-2k)+\frac{1}{2}(n-2k)^2-\frac{1}{2}n=\frac{1}{2}(n-2k)^2+\frac{1}{2}(n-2k)-k$. This is increasing as a function of $n-2k$, and so the minimum in eqn.~\eqref{K_n_D1,2basic} is achieved when $n-2k=\lc \sqrt{n} \rc$ or $\lc \sqrt{n} \rc+1$, depending on which of these has the same parity as $n$. Hence, the minimum is achieved when $k=\frac{1}{2}\left(n-\left\lc \sqrt{n} \right\rc\right)$ or $\frac{1}{2}\left(n-1-\left\lc \sqrt{n} \right\rc\right)$, whichever is an integer. 
\vskip 5pt

\noindent {\sc Case} II. If $n-2k<\sqrt{n}$, then the expression in eqn.~\eqref{K_n_D1,2basic} is $\frac{1}{2}n-\frac{1}{2}(n-2k)^2+(n-2k)=\frac{1}{2}(n+1)-\frac{1}{2}\left( (n-2k)^2 -2(n-2k)+1 \right)=\frac{1}{2}(n+1)-\frac{1}{2}(n-2k-1)^2$. This is decreasing as a function of $n-2k$, and so the minimum in eqn.~\eqref{K_n_D1,2basic} is achieved when $n-2k=\lf \sqrt{n} \rf$ or $\lf \sqrt{n} \rf-1$, depending on which of these has the same parity as $n$, if $n$ is not a square. If $n$ is a square, then the minimum is achieved when $n-2k=\sqrt{n}-2$, which has the same parity as $n$. Hence, the minimum is achieved when $k=\frac{1}{2}\left(n-\left\lf \sqrt{n} \right\rf\right)$ or $\frac{1}{2}\left(n+1-\left\lf \sqrt{n} \right\rf\right)$, whichever is an integer, when $n$ is not a square, and when $k=\frac{1}{2}(n+2-\sqrt{n})$ when $n$ is a square. 
\vskip 5pt

The minimum in eqn.~\eqref{K_n_D1,2basic} is the smaller of the two minima we considered in {\sc Cases} I and II.  
\vskip 5pt
 
If $n=a^2$, then $\lf \sqrt{n} \rf=\lc \sqrt{n} \rc=a$. Since $a,n$ have the same parity, the minimum in eqn.~\eqref{K_n_D1,2basic} is achieved when $n-2k$ equals either $a$ ({\sc Case} I) or $a-2$ ({\sc Case} II), by the previous arguments. The corresponding values of $(n-2k)+\frac{1}{2}\left| (n-2k)^2-n \right|$ are $a$ and $a-2+(2a-2)=3a-4$, and the smaller of these equals $a=\sqrt{n}$.  

Now suppose $a^2<n<(a+1)^2$. If $a,n$ have the same parity, the minimum in eqn.~\eqref{K_n_D1,2basic} is achieved when $n-2k$ equals either $a+2$ ({\sc Case} I) or $a$ ({\sc Case} II), by the previous arguments. The corresponding values of $(n-2k)+\frac{1}{2}\left| (n-2k)^2-n \right|$ are $a+2+\frac{1}{2}\big((a+2)^2-n\big)$ and $a+\frac{1}{2}(n-a^2)$, and the smaller of these equals $a+\frac{1}{2}(n-a^2)$. 

If $a,n$ have opposite parity, the minimum in eqn.~\eqref{K_n_D1,2basic} is achieved when $n-2k$ equals either $a+1$ ({\sc Case} I) or $a-1$ ({\sc Case} II). The corresponding values of $(n-2k)+\frac{1}{2}\left| (n-2k)^2-n \right|$ are $a+1+\frac{1}{2}\big((a+1)^2-n\big)$ and $a-1+\frac{1}{2}\big(n-(a-1)^2\big)$; these are equal if $n=a^2+3$, the smaller of these is $a-1+\frac{1}{2}\big(n-(a-1)^2\big)$ if $n=a^2+1$, and $a+1+\frac{1}{2}\big((a+1)^2-n\big)$ if $n \in \{a^2+5,a^2+7,a^2+9,\ldots,a^2+2a-1\}$. 

Therefore, 
\begin{equation*} 
{\D}_1(G) = \begin{cases} 
a+\frac{1}{2}(n-a^2) & \:\:\mbox{if}\:\: n \in \{a^2+2t: 0 \le t \le a\}, \\
2a-1 & \:\:\mbox{if}\:\: n=a^2+1, \\
a+1+\frac{1}{2}\big((a+1)^2-n\big) & \:\:\mbox{if}\:\: n \in \{a^2+2t+1: 1 \le t \le a-1\}. \\
\end{cases} 
\end{equation*}
\end{Pf}
\vskip 10pt

\section{The complete $r$-partite graphs ${\K}_{n_1,\ldots,n_r}$} \label{r-partite}
\vskip 10pt

Cahit \cite{Cah87} proved that the complete bipartite graph ${\K}_{m,n}$ is cordial for every pair of positive integers $m,n$. We find upper bounds for ${\D}_1$ for all complete $r$-partite graphs ${\K}_{n_1,\ldots,n_r}$ and exactly determine ${\D}_2$ in all cases. 
\vskip 5pt

\begin{thm} \label{r_partite}
If $n_1,\ldots,n_r$ are positive integers of which $s$ are odd and $(2a)^2 \le s < (2a+2)^2$, then
\begin{itemize}
\item[{\rm (i)}]  
\[ {\D}_1\big({\K}_{n_1,\ldots,n_r}\big) \begin{cases}
                                                              = \sqrt{s} & \text{ if } s \text{ is a square}, \\
                                                              \le \frac{1}{2} \big(s+1 - (2a-1)^2\big) & \text{ if } (2a)^2 < s < (2a+1)^2 \text{ and } s \text{ is even}, \\
                                                              \le \frac{1}{2} \big(s+1 - (2a-2)^2\big) & \text{ if } (2a)^2 < s < (2a+1)^2 \text{ and } s \text{ is odd}, \\
                                                              \le \frac{1}{2} \big(s+1 - (2a-1)^2\big) & \text{ if } (2a+1)^2 < s < (2a+2)^2 \text{ and } s \text{ is even}, \\
                                                              \le \frac{1}{2} \big(s+1 - (2a)^2\big) & \text{ if } (2a+1)^2 < s < (2a+2)^2 \text{ and } s \text{ is odd}. 
                                                            \end{cases}
\]
Moreover, we have 
\[ \lf \sqrt{s} \rf \le {\D}_1\big({\K}_{n_1,\ldots,n_r}\big) \le 3 \lf \sqrt{s} \rf. \]
\item[{\rm (ii)}]
\[ {\D}_2\big({\K}_{n_1,\ldots,n_r}\big) = \left\lf \tfrac{s}{2} \right\rf. \]
\end{itemize}
\end{thm}

\begin{Pf}
Let $n_1,\ldots,n_r$ be positive integers, of which $n_1,\ldots,n_s$ are odd. Let $G={\K}_{n_1,\ldots,n_r}$, with partite sets $X_i$, $|X_i|=n_i$, $1 \le i \le r$. Let $f: V(G) \to \{0,1\}$, and let $k_i$ denote the number of vertices labelled $0$ in $X_i$, $ 1 \le i \le r$. Thus,  
\begin{eqnarray*} \label{r_partite_del_v}
v_0(f) = \sum_{i=1}^r k_i, \;\; v_1(f) = \sum_{i=1}^r (n_i-k_i); \quad {\De}_v(f) = \left| \sum_{i=1}^r (n_i-2k_i) \right|.
\end{eqnarray*}
\begin{eqnarray*} \label{r_partite_e0e1}
e_0(f) = \sum_{1 \le i<j \le r} \big( k_ik_j + (n_i-k_i)(n_j-k_j) \big), \;\; e_1(f) = \sum_{1 \le i<j \le r} \big( k_i(n_j-k_j) + (n_i-k_i)k_j \big);
\end{eqnarray*}
\begin{eqnarray*} \label{r_partite_del_e}
{\De}_e(f) = \left| \sum_{1 \le i<j \le r} \big( n_in_j -2(n_ik_j+n_jk_i) + 4k_ik_j \big) \right| = \left| \sum_{1 \le i<j \le r} (n_i-2k_i)(n_j-2k_j) \right|. 
\end{eqnarray*}
\vskip 5pt

\noindent Write $d_i=n_i-2k_i$, $1 \le i \le r$. Then  
\begin{eqnarray} \label{r_partite_D1,2basic}
{\D}_1\big({\K}_{n_1,\ldots,n_r}\big) & = & \min_{\substack{0 \le k_i \le n_i \\ 1 \le i \le r}} \left( \left| \sum_{i=1}^r d_i \right| + \left| \sum_{1 \le i<j \le r} d_i d_j \right| \right) \nonumber \\
& = & \min_{\substack{0 \le k_i \le n_i \\ 1 \le i \le r}} \left( \left| \sum_{i=1}^r d_i \right| + \frac{1}{2} \:\: \left| \left( \sum_{i=1}^r d_i \right)^2 - \sum_{i=1}^r d_i^2 \right| \right). \\ \label{r_partite_D1}
{\D}_2\big({\K}_{n_1,\ldots,n_r}\big) & = & \min_{\substack{0 \le k_i \le n_i \\ 1 \le i \le r}} \:\: \left\{  \frac{1}{2} \left| \left( \sum_{i=1}^r d_i \right)^2 - \sum_{i=1}^r d_i^2 \right|: \left| \sum_{i=1}^r d_i \right| \le 1 \right\}. \label{r_partite_D2}
\end{eqnarray}

\begin{itemize}
\item[{\rm (i)}]
We provide upper bounds for $\D_1({\K}_{n_1, \ldots, n_r})$ by taking the two cases, $s$ is even and $s$ is odd, and using eqn.~\eqref{r_partite_D1}. Observe that $\sum_{i=1}^rd_i^2 \ge s$ since $|d_i| \ge 1$ if $n_1$ is odd. Thus, if $|\sum_{i=1}^rd_i| \le \lf \sqrt{s} \rf$, then
\begin{eqnarray}
\left| \sum_{i=1}^r d_i \right| + \frac{1}{2} \:\: \left| \left( \sum_{i=1}^r d_i \right)^2 - \sum_{i=1}^r d_i^2 \right| & = & \left| \sum_{i=1}^r d_i \right| - \frac{1}{2}\left| \sum_{i=1}^r d_i \right|^2 + \frac{1}{2}\sum_{i=1}^r d_i^2 \nonumber \\
& = & \frac{1}{2}\left( \sum_{i=1}^r d_i^2 + 1 - \left( \left| \sum_{i=1}^r d_i \right| - 1 \right)^2 \right). \label{eqn:d_1_expr_case}
\end{eqnarray}
Consider a labelling $f$ such that
\[ d_i = \begin{cases}
              -1 & \text{ if } i \le \left\lc \frac{s - \lf \sqrt{s} \rf}{2} \right\rc, \\
               1 & \text{ if } \left\lc \frac{s - \lf \sqrt{s} \rf}{2} \right\rc < i \le s, \\
               0 & \text{ if } s < i \le r.
             \end{cases}
\]
Notice that for this labelling, 
\begin{eqnarray}
\sum_{i=1}^r d_i^2 = s \quad \text{and} \quad \sum_{i=1}^rd_i = s - 2\left\lc \frac{s - \lf \sqrt{s} \rf}{2} \right\rc \le \lf \sqrt{s} \rf. \label{eqn:values_for_f}
\end{eqnarray}
Eqns.~\eqref{eqn:d_1_expr_case} and \eqref{eqn:values_for_f} together give the upper bounds for $\D_1\big({\K}_{n_1, \ldots, n_r}\big)$. When $s=t^2$, 
\[ \sum_{i=1}^r d_i = s - 2 \left\lc \frac{s - \lf \sqrt{s} \rf}{2} \right\rc = t^2 - 2 \left\lc \frac{t^2-t}{2} \right\rc = t^2 - 2 \left( \frac{t^2-t}{2} \right) = t. \] 
Substituting $\sum_{i=1}^r d_i^2 = t^2$ in eqn.~\eqref{eqn:d_1_expr_case} gives the upper bound as $\frac{1}{2}\left(t^2+1-(t-1)^2\right)=t$. To show that $t$ is also a lower bound, we consider two cases given later. 
\vskip 5pt

\noindent Now consider the four cases when $s$ is not a square. Substituting $\sum_{i=1}^r d_i^2 = s$ in eqn.~\eqref{eqn:d_1_expr_case} and comparing with the upper bounds for $\D_1({\K}_{n_1, \ldots, n_r})$ to be proven, we need to show that $\sum_{i=1}^r d_i = s-2\left\lc \frac{s - \lf \sqrt{s} \rf}{2} \right\rc$ equals $2a, 2a-1, 2a, 2a+1$ in these four cases. We show the first such case; the other three cases are similarly derived. For the case $(2a)^2<s<(2a+1)^2$, $s$ even, we have 
\[ \sum_{i=1}^r d_i = s - 2 \left\lc \frac{s - \lf \sqrt{s} \rf}{2} \right\rc = s - 2 \left\lc \frac{s-2a}{2} \right\rc = s - 2\left( \frac{s-2a}{2} \right) = 2a. \]
\vskip 5pt

\noindent This proves the upper bound for $\D_1({\K}_{n_1, \ldots, n_r})$ in all cases. Further, when $s$ is not a square, we can combine the four upper bounds to get a common upper bound. We note that in the third case,  
\[ \frac{1}{2} \left( s+1 - (2a-1)^2 \right) \le \frac{1}{2} \left( (2a+2)^2 - (2a-1)^2 \right) = \frac{3}{2}(4a+1) < 3 \left\lf \sqrt{s} \right\rf. \]
A similar computation shows the upper bounds in each of the other cases to be no more than $3 \lf \sqrt{s} \rf$. 
\vskip 5pt

For the lower bound for $\D_1({\K}_{n_1, \ldots, n_r})$, we show that $\D_1({\K}_{n_1, \ldots, n_r}) \ge \lf \sqrt{s} \rf$. We consider two cases. 
\vskip 5pt

\noindent {\sc Case} I. If $|\sum_{i=1}^r d_i| \ge \lf \sqrt{s} \rf$, then
\[ \left| \sum_{i=1}^r d_i \right| + \frac{1}{2} \:\: \left| \left( \sum_{i=1}^r d_i \right)^2 - \sum_{i=1}^r d_i^2 \right| \ge \lf \sqrt{s} \rf. \]
\vskip 5pt

\noindent {\sc Case} II. If $|\sum_{i=1}^r d_i| < \lf \sqrt{s} \rf$, then by eqn.~\eqref{eqn:d_1_expr_case} and using the fact that $\sum_{i=1}^r d_i^2 \ge s$,
\begin{eqnarray}
\left| \sum_{i=1}^r d_i \right| + \frac{1}{2} \:\: \left| \left( \sum_{i=1}^r d_i \right)^2 - \sum_{i=1}^r d_i^2 \right| & = & \frac{1}{2}\left( \sum_{i=1}^r d_i^2 + 1 - \left( \left| \sum_{i=1}^r d_i \right| - 1 \right)^2 \right) \nonumber \\
& \ge & \frac{1}{2} \left(s + 1 - \left( \left| \sum_{i=1}^r d_i \right| - 1 \right)^2 \right) \nonumber \\
& > & \frac{1}{2} \left(s + 1 - \left( \lf \sqrt{s} \rf - 1 \right)^2 \right) \nonumber \\
& \ge & \lf \sqrt{s} \rf. 
\end{eqnarray}
\vskip 5pt

\noindent From {\sc Case} I and II, we see that $\D_1({\K}_{n_1, \ldots, n_r}) \ge \lf \sqrt{s} \rf$. In particular, this proves $\D_1({\K}_{n_1, \ldots, n_r})=\sqrt{s}$ when $s$ is a square. 

\item[{\rm (ii)}]
We determine $\D_2({\K}_{n_1, \ldots, n_r})$ by taking the two cases $s$ is even and $s$ is odd and using eqn.~\eqref{r_partite_D2}. 
\vskip 5pt

\noindent {\sc Case} I. ($s=2q$) In this case, since the number of vertices in $G$ is even, any labelling $f$ of $\D_2\big({\K}_{n_1, \ldots, n_r}\big)$ must satisfy ${\De}_v(f) = 0$. This implies that
\begin{eqnarray} \label{r_partite_D2lower_even}
\D_2({\K}_{n_1, \ldots, n_r}) = \min_{\substack{0 \le k_i \le n_i \\ 1 \le i \le r}} \left\{ \frac{1}{2} \sum_{i=1}^r d_i^2: \sum_{i=1}^r d_i = 0\right\} \ge \frac{1}{2} \left( \sum_{i=1}^r\min_{0 \le k_i \le n_i} d_i^2 \right) \ge \frac{s}{2}. 
\end{eqnarray}
To show that the lower bound in eqn.~\eqref{r_partite_D2lower_even} can be achieved, consider $g: V(G) \to \{0, 1\}$ for which  
\[ k_i = \begin{cases} 
              \left\lf \frac{n_i}{2} \right\rf & \:\mbox{ if }\: i \in \{1, \ldots, q\} \bigcup \{s+1, \ldots, r\}, \\
              \left\lf \frac{n_i}{2} \right\rf + 1 & \:\mbox{ if }\: i \in \{q+1, \ldots, s\}.
             \end{cases}
\]
Then ${\De}_v(g) = 0$ and ${\De}_e(g) = s/2$, so that $\D_2\big({\K}_{n_1, \ldots, n_r}\big) = s/2$ for this case.
\vskip 5pt

\noindent {\sc Case} II. ($s = 2q+1$) In this case, since the number of vertices in $G$ is odd, any labelling $f$ of $\D_2\big({\K}_{n_1, \ldots, n_r}\big)$ must satisfy ${\De}_v(f) = 1$. Again, since $s$ is odd, $\sum_{i=1}^r d_i^2 \ge 1$, so that 
\begin{eqnarray} \label{r_partite_D2lower_odd}
\D_2\big({\K}_{n_1, \ldots, n_r}\big) = \min_{\substack{0 \le k_i \le n_i \\ 1 \le i \le r}} \left\{ \frac{1}{2} \left( \sum_{i=1}^r d_i^2 - 1 \right): \sum_{i=1}^r d_i = \pm 1 \right\} \ge \frac{1}{2} \left( \sum_{i=1}^r \min_{0 \le k_i \le n_i} d_i^2 - 1 \right) \ge \frac{s-1}{2}.
\end{eqnarray}
To show that the lower bound in eqn.~\eqref{r_partite_D2lower_odd} can be achieved, note that the labelling in {\sc Case} I satisfies ${\De}_v(g) = 1$ and ${\De}_e(g) = (s-1)/2$, so that $\D_2\big({\K}_{n_1, \ldots, n_r}\big) = (s-1)/2$ for this case.
\end{itemize}
\end{Pf}
\vskip 5pt

\begin{cor} {\bf (Lee \& Liu \cite{LL91})} \\[5pt]
A complete multipartite graph is cordial if and only if at most three of its partite sets have odd cardinality. 
\end{cor}

\begin{Pf}
This follows immediately from Theorem \ref{r_partite}, part (ii) and the fact that $G$ is cordial if and only if $\D_2(G) \le 1$.  
\end{Pf}
\vskip 10pt

\section{The cycles ${\C}_n$} \label{cycle}
\vskip 10pt

Cahit \cite{Cah87} proved that an Eulerian graph is not cordial if its size is congruent to $2$ modulo $4$. In particular, the cycle ${\C}_n$ is not cordial if $n \equiv 2\pmod{4}$. We determine ${\D}_1({\C}_n)$ and ${\D}_2({\C}_n)$ for each positive integer $n$. 
\vskip 5pt

\begin{thm} \label{C_n}
\[ {\D}_1\big({\C}_n\big) = \begin{cases} 
                                             0 & \mbox{ if } n \equiv 0\!\!\!\!\pmod{4}, \\
                                             2 & \mbox{ if } n \not\equiv 0\!\!\!\!\pmod{4}, 
                                            \end{cases}
\quad {\rm and} \quad 
{\D}_2\big({\C}_n\big) = \begin{cases} 
                                          0 & \mbox{ if } n \equiv 0\!\!\!\!\pmod{4}, \\
                                          1 & \mbox{ if } n \equiv 1,3\!\!\!\!\pmod{4}, \\
                                          2 & \mbox{ if } n \equiv 2\!\!\!\!\pmod{4}. 
                                        \end{cases}
\] 
\end{thm}

\begin{Pf}
Let the vertices of ${\C}_n$ be $v_1,\ldots,v_n$ and the edges $e_1,\ldots,e_n$, where $e_i=v_iv_{i+1}$ for $1 \le i \le n-1$ and $e_n=v_nv_1$. Assign $0$ to $v_1$ and $v_2$, so that $0$ is also assigned to $e_1$. If labels have been assigned to $v_1,\ldots,v_k$, $1 \le k \le n-1$, and hence to $e_1,\ldots,e_{k-1}$, assign $0$ or $1$ to $v_{k+1}$ such that the induced labelling on $e_k=v_kv_{k+1}$ is not the label assigned to $e_{k-1}$. So the sequence  
\[ 0, 0, 1, 1, 0, 0, 1, 1, 0, 0, 1, 1, 0, 0, 1, 1, \ldots \]
of labels assigned to vertices induces an alternating sequence of $0$'s and $1$'s to the edges. This labelling may be given by
\begin{eqnarray} \label{f_cycle1}
f(v_i) = \begin{cases} 
              0 \:\:\mbox{if}\:\: i \equiv \:1,2\!\!\pmod{4}, \\ 
              1 \:\:\mbox{if}\:\: i \equiv \:0,3\!\!\pmod{4}. 
             \end{cases} 
\end{eqnarray}  
Thus, $v_0(f)-v_1(f)$ equals $0,1,2,1$, and $e_0(f)-e_1(f)$ equals $0,1,2,-1$ for $n=4k,4k+1,4k+2,4k+3$, respectively. We consider the cases $n=4k+2$ and $n \ne 4k+2$ separately. 
\vskip 5pt

\noindent {\sc Case} I. ($n \ne 4k+2$) The above example shows ${\D}_1\big({\C}_{4k}\big)={\D}_2\big({\C}_{4k}\big)=0$. If $n$ is odd, then the above example shows ${\De}_v(f)={\De}_e(f)=1$ since the number of vertices as well as the number of edges labelled $0$ and $1$ cannot be equal. This proves the theorem for these cases. 
\vskip 5pt

\noindent {\sc Case} II. ($n=4k+2$) Let $f: \{v_1,\ldots,v_{4k+2}\} \to \{0,1\}$ be any labelling, and let $\ov{f}: \{e_1,\ldots,e_{4k+2}\} \to \{0,1\}$ be the labelling on the edges induced by $f$. For any cycle $v_1 \dots v_n$ where $v_1 = v_n$,
\begin{eqnarray} \label{eqn:e1_cycle_parity}
e_1(f) = \sum_{i=1}^{n} \ov{f}(e_i) = \sum_{i=1}^{n} \big| f(v_i)-f(v_{i+1}) \big| \equiv \sum_{i=1}^{n} \big( f(v_i)-f(v_{i+1}) \big) = 0 \pmod{2}.
\end{eqnarray}
This shows ${\D}_2\big({\C}_{4k+2}\big) \ge {\D}_1\big({\C}_{4k+2}\big) \ge {\De}_e(f) \ge 2$. 
\vskip 5pt

The labelling $f$ defined above satisfies $v_0(f)-v_1(f)=2$ for $n=4k+2$. For this labelling, $f(v_{4k+2})=0$. If we instead define $f(v_{4k+2})=1$, then $v_0(f)$ decreases by one and $v_1(f)$ increases by one. Thus, the modified function $\t{f}$ satisfies $v_0(\t{f})=v_1(\t{f})$. The changes in assignment of labels to edges is due to the sequence of labels for $v_{4k+1},v_{4k+2},v_1$ changing from $0,0,1$ to $0,1,1$. This results in no change in $e_0(f)$ or $e_1(f)$, so $e_0(\t{f})-e_1(\t{f})=2$. Therefore, $\t{f}$ is a labelling satisfying ${\De}_v(\t{f})=0$ and ${\De}_e(\t{f})=2$. This shows ${\D}_2\big({\C}_{4k+2}\big)={\D}_1\big({\C}_{4k+2}\big)=2$, completing the proof of the theorem. 
\end{Pf}
\vskip 5pt

\begin{cor} \label{eulerian} {\bf (Cahit \cite{Cah87})} \\[5pt]
An Eulerian graph $G$ with $4k+2$ edges is not cordial for all positive integer $k$.
\end{cor}

\begin{Pf}
Since $G$ has $4k+2$ edges, if $f: V(G) \to \{0,1\}$ is a cordial labelling then $e_0(f) = e_1(f) = 2k+1$. Let $v_1,\ldots,v_{4k+3}$ be an Eulerian circuit of $G$ such that $v_{4k+3} = v$. Then, by eqn. \eqref{eqn:e1_cycle_parity}, $e_1(f)$ has to be even for every labelling $g: V(G) \to \{0,1\}$. Thus, no such cordial labelling exists.
\end{Pf}
\vskip 10pt

\section{The wheel graphs ${\W}_n$} \label{wheel}
\vskip 10pt

The wheel graph ${\W}_n$ is the $n$-vertex graph with $n-1$ vertices forming a cycle ${\C}_{n-1}$ and a central vertex $x$ adjacent to every vertex on the cycle. Thus, ${\W}_n$ has $2n-2$ edges. Cahit \cite{Cah87} proved that ${\W}_n$ is cordial if and only if $n \not \equiv 3\pmod{4}$. We determine ${\D}_1({\W}_n)$ and ${\D}_2({\W}_n)$ for each positive integer $n$. 
\vskip 5pt

\begin{thm} \label{W_n}
\[ {\D}_1\big({\W}_n\big) = \begin{cases} 
                                             0 & \mbox{if}\:\: n \equiv 2\!\!\!\!\pmod{4}, \\
                                             1 & \mbox{if}\:\: n \equiv 1,3\!\!\!\!\pmod{4}, \\
                                             2 & \mbox{if}\:\: n \equiv 0\!\!\!\!\pmod{4}, 
                                           \end{cases}
\quad {\rm and} \quad 
{\D}_2\big({\W}_n\big) = \begin{cases} 
                                          0 & \mbox{if}\:\: n \not\equiv 0\!\!\!\!\pmod{4}, \\
                                          2 & \mbox{if}\:\: n \equiv 0\!\!\!\!\pmod{4}. 
                                        \end{cases}
\] 
\end{thm}

\begin{Pf}
Consider the function $f$ defined by eqn.~\eqref{f_cycle1} in Theorem \ref{C_n}. We label ${\W}_n$ as follows. 
\begin{eqnarray*} \label{f_cycle2}
g(v) = \begin{cases} 
              f(v) \quad\mbox{if}\:\: v \in {\C}_{n-1}, \\ 
              1 \quad\mbox{if}\:\: v = x. 
             \end{cases} 
\end{eqnarray*}  
The values of $v_0(f) - v_1(f)$ and $e_0(f) - e_1(f)$ for ${\C}_{n-1}$ are tabulated in Table \ref{wheels}. The structure of the wheel and definition of $g$ imply the following.
\begin{eqnarray*}
v_0(g) - v_1(g) & = & v_0(f) - v_1(f) - 1, \\
e_0(g) - e_1(g) & = & \big( e_0(f) - e_1(f) \big) - \big( v_0(f) - v_1(f) \big).
\end{eqnarray*}
These values are tabulated in Table \ref{wheels}. The ``minimality'' of the labelling $g$ when $n \ne 4k$ is implied by the parity of edge counts and vertex counts. This proves the theorem except when $n = 4k$. 

\begin{table}[H] 
\centering
\renewcommand{\arraystretch}{0.9}
\begin{tabular}{|c|c|c|c|c|c|c|} \hline
{\sc Subcase} & $v_0(f) - v_1(f)$ & $e_0(f) - e_1(f)$ & $v_0(g) - v_1(g)$ & $e_0(g) - e_1(g)$ & ${\D}_1\big({\W}_n\big)$ & ${\D}_2\big({\W}_n\big)$ \\ \hline
$n = 4k$ & $1$ & $-1$ & $0$ & $-2$ & $2$ & $2$ \\ \hline
$n = 4k+1$ & $0$ & $0$ & $-1$ & $0$ & $1$ & $0$ \\ \hline
$n = 4k+2$ & $1$ & $1$ & $0$ & $0$ & $0$ & $0$ \\ \hline
$n = 4k+3$ & $2$ & $2$ & $1$ & $0$ & $1$ & $0$ \\ \hline
\end{tabular}
\caption{Calculation of ${\D}_1\big({\W}_n\big)$ and ${\D}_2\big({\W}_n\big)$ using the function $g$} \label{wheels}
\end{table}
\vskip 5pt

\noindent The minimality for the case when $n=4k$ can be argued as follows. Suppose there exists a labelling $\tilde{g}$ of ${\W}_{4k}$ such that $v_0(\tilde{g}) = v_1(\tilde{g}) = 2k$ and $e_0(\tilde{g}) = e_1(\tilde{g}) = 4k-1$. Without loss of generality, let $\tilde{g}(x) = 1$. If $\tilde{f}$ is the labelling induced by $\tilde{g}$ on ${\C}_{4k-1}$ then $v_0(\tilde{f}) = 2k$ and hence $e_1(\tilde{f}) = e_1(\tilde{g}) - v_0(\tilde{f}) = 2k-1$, which contradicts eqn. \eqref{eqn:e1_cycle_parity}.  Thus, there exists no such labelling and by parity arguments, $g$ is a ``minimal'' labelling. 
\end{Pf}
\vskip 10pt

\section{The Fan graphs ${\F}_{m, n}$} \label{fan}
\vskip 10pt

The fan graph $\F_{m,n}$ is the join between a path $\P_n$ and an empty graph $\overline{\K}_m$. Thus, it has $m+n$ vertices and $mn+n-1$ edges. Cahit \cite{Cah87} proved that all fans are cordial, so that ${\D}_1({\F}_{m, n}) \le 2$ and ${\D}_2({\F}_{m, n}) \le 1$. We exploit the proof of cordiality to compute the exact values.  

\begin{thm} \label{F_m,n}
\begin{eqnarray*}
 {\D}_1\big({\F}_{m,n}\big) = \begin{cases} 
                                               2 & \mbox{ if $m$ is odd, $n$ is even}, \\
                                               1 & \mbox{ otherwise,}
                                              \end{cases} 
\qquad \text{and} \qquad 
{\D}_2\big({\F}_{m,n}\big) = \begin{cases} 
                                              0 & \mbox{ if $m$ is even, $n$ is odd}, \\
                                              1 & \mbox{ otherwise.}
                                             \end{cases}
\end{eqnarray*}
\end{thm}

\begin{Pf}
Let the fan graph $\F_{m,n}$ denote the join of $\P_n = v_1\ldots v_n$ and $\overline{\K}_m = \{u_1, \ldots, u_m\}$. Hence $E(\F_{m, n}) = E(\P_n) \cup (V(\P_n) \times V(\ov{\K}_m))$. We define the following labelling $f: V(G) \to \{0, 1\}$.
\begin{eqnarray}
& f(v_i) = \begin{cases}
                  0 & \mbox{ if } \:\: i \equiv 0,1 \!\!\!\!\pmod{4}, \\
                  1 & \mbox{ if } \:\: i \equiv 2,3 \!\!\!\!\pmod{4},
                \end{cases}
\quad {\rm and} \quad 
f(u_i) = \begin{cases}
               0 & \mbox{ if } \:\: i \in \{1, \ldots \lf m/2 \rf\}, \\
               1 & \mbox{ if } \:\: i \in \{\lc m/2 \rc + 1, \ldots, m\},
             \end{cases} \\
& f(u_{\lc m/2 \rc}) = \begin{cases}
                                     0 & \mbox{ if } \:\: m \equiv 1 \!\!\!\!\pmod{2} \text{ and } n \equiv 0, 3 \!\!\!\!\!\pmod{4}, \\
                                     1 & \mbox{ if } \:\: m \equiv 1 \!\!\!\!\pmod{2} \text{ and } n \equiv 1, 2 \!\!\!\!\pmod{4}.
                                   \end{cases}
\end{eqnarray}
\vskip 5pt
We now compute the values of $v_0(f) - v_1(f)$ and $e_0(f) - e_1(f)$ to show that \(f\) yields the stated \({\D}_1\) and \({\D}_2\) values. The minimality follows by parity.
\vskip 5pt

\noindent {\sc Case} I. ($m$ is even) The labelling $f$ labels equal number of vertices in ${\ov{\K}}_m$ with $0$ and $1$. Thus, there are an equal number of $0$ and $1$ labelled edges in $V(\P_n) \times V(\ov{\K}_m)$. Hence, $v_0(f) - v_1(f) = v_0(f|_{\P_n}) - v_1(f|_{\P_n})$ and $e_0(f) - e_1(f) = e_0(f|_{\P_n}) - e_1(f|_{\P_n})$ where $f|_{G}$ is the labelling $f$ restricted to a subgraph $G$. An easy induction on $n$ yields the following.
\begin{eqnarray}
v_0(f|_{\P_n}) - v_1(f|_{\P_n}) = \begin{cases}
                                                  0 & \text{ if } n \equiv 0 \!\!\!\!\pmod{2}, \\
                                                  1 & \text{ if } n \equiv 1 \!\!\!\!\pmod{4}, \\
                                                -1 & \text{ if } n \equiv 3 \!\!\!\!\pmod{4}, 
                                                 \end{cases} \\
e_0(f|_{\P_n}) - e_1(f|_{\P_n}) = \begin{cases}
                                                   0 & \text{ if } n \equiv 1 \!\!\!\!\pmod{2}, \\
                                                  -1 & \text{ if } n \equiv 0 \!\!\!\!\pmod{2}. 
                                                 \end{cases}
\end{eqnarray}
\vskip 5pt

\noindent {\sc Case} II. ($m$ is odd, $n \equiv 1 \pmod{4}$) Removing $u_{\lc m/2 \rc}$ reduces the problem to an instance of {\sc Case} I. Adding back $u_{\lc m/2 \rc}$ gives us $v_0(f) - v_1(f) = v_0(f|_{\F_{m-1,n}}) - v_1(f|_{\F_{m-1,n}}) - 1$ and $e_0(f) - e_1(f) = e_0(f|_{\F_{m-1,n}}) - e_1(f|_{\F_{m-1,n}}) - 1$. Thus, $v_0(f) - v_1(f) = 0$ and $e_0(f) - e_1(f) = -1$.
\vskip 5pt

\noindent {\sc Case} III. ($m$ is odd, $n \equiv 3 \pmod{4}$) Removing $u_{\lc m/2 \rc}$ reduces the problem to an instance of {\sc Case} I. Adding back $u_{\lc m/2 \rc}$ gives us $v_0(f) - v_1(f) = v_0(f|_{\F_{m-1,n}}) - v_1(f|_{\F_{m-1,n}}) + 1$ and $e_0(f) - e_1(f) = e_0(f|_{\F_{m-1,n}}) - e_1(f|_{\F_{m-1,n}}) - 1$. Thus, $v_0(f) - v_1(f) = 0$ and $e_0(f) - e_1(f) = -1$. 
\end{Pf}
\vskip 10pt

\section{Concluding Remarks} \label{CR}
\vskip 10pt

We have introduced two measures of cordiality, and have investigated their values for several important classes of graphs. The interested reader is invited to extend these findings to a larger collection of graphs. We close this paper by listing three directions of further research. 
\vskip 5pt

\noindent {\bf Open Problem 1.} Determine the exact value of $\D_1({\K}_{n_1, \ldots, n_r})$ in those cases where the number of odd sized parts is not a perfect square. 
\vskip 5pt

\noindent {\bf Open Problem 2.} Among all graphs $G$ of order $n$, determine the largest possible value of ${\D}_i(G)$, with $i=1,2$. Also, find all such extremal graphs in both cases. 
\vskip 5pt

\noindent {\bf Open Problem 3.} Given $n$ and $D_1$, determine the maximum size of a graph $G$ with $|V(G)|=n$ and ${\D}_1(G)=D_1$. The same problem for ${\D}_2$. Also, find all such extremal graphs in both cases. 
\vskip 20pt

\noindent {\bf Acknowledgement.} The authors gratefully acknowledge the comments of the two reviewers.

\end{document}